%%%%%%%%%%%%%%%%%%%%%%%%%%%%%%%%%
%
%     
%
% Version of   March 29,2008
%
% Last changes by DR
%
%%%%%%%%%%%%%%%%%%%%%%%%%%%%%%%%%
\input amstex.tex
\input amsppt.sty
%%%%%%%%%%%%%%%%%%%%%%%%%%%%%%
\input pictex
%%%%%%%%%%%%%%%%%%%%%%%%%%%%%%
\documentstyle{amsppt}
%\magnification=\magstep{1.5} %=1200
%\pagewidth{6.5truein} \pageheight{9.5truein} \TagsOnRight
\NoBlackBoxes
\hoffset=5mm \voffset-10mm
\let\e\varepsilon
\let\phi\varphi
\define\N{{\Bbb N}}
\define\R{{\Bbb R}}
\define\FF{{\Cal F}}
\define\CC{{\Cal C}}

\def\ss{\subset}

\let\i\infty
\def\d{\delta}

\def\l{\lambda}
\def\a{\alpha}
\def\b{\beta}
\def\q{\qed}
\def\f{\frac}

\def\conv{\hbox{conv}}
\def\dist{\hbox{dist}}
\def\aff{\hbox{aff}}
\def\Ker{\hbox{Ker}}
\def\Gr{\hbox{Gr}}
\def\int{\hbox{int}}
\def\Im{\hbox{Im}}
\def\span{\hbox{span}}

\topmatter
\title  Hereditary invertible  linear surjections and splitting problems for selections \endtitle

\author       Du\v{s}an Repov\v{s} and Pavel V. Semenov   \endauthor
\leftheadtext{Du\v{s}an Repov\v{s} and Pavel V. Semenov}
\rightheadtext{Hereditary invertible  linear surjections...}

\address
Institute of Mathematics, Physics and Mechanics, and Faculty of Education, University of
Ljubljana, P. O. Box 2964, Ljubljana, Slovenia 1001
\endaddress
\email dusan.repovs\@guest.arnes.si
\endemail
\address
Department of Mathematics, Moscow City Pedagogical University,
2-nd Selsko\-khozyast\-vennyi pr.\,4,\,Moscow,\,Russia 129226
\endaddress
\email pavels\@orc.ru
\endemail

\abstract  Let $A+B$ be the pointwise (Minkowski) sum of two
convex subsets $A$ and $B$ of a Banach space. Is it true that
every continuous mapping $h:X \to A+B$ splits into a sum $h=f+g$
of continuous mappings $f:X \to A$ and $g:X \to B$? We study this
question within a wider framework of splitting techniques of
continuous selections. Existence of splittings is guaranteed by
hereditary invertibility of linear surjections between Banach
spaces. Some affirmative and negative results on such
invertibility with respect to an appropriate class of convex compacta are
presented. As a corollary, a positive answer to the above question
is obtained for strictly convex finite-dimensional precompact
spaces.
\endabstract
\endtopmatter

\document
%\baselineskip=6truemm

\head {\bf 1. Introduction}\endhead

Recall that a single-valued mapping $f:X \to Y$ is said to be a
{\it selection} of a multivalued mapping $F: X \to Y$ provided
that $f(x) \in F(x),$  for every $x \in X$. Classically,
selections exist in the category of topological spaces (for
details see \cite{M1, Mi, RS}), or in the category of measurable
spaces (see \cite{AC, AF, RS}). Here we shall restrict ourselves
only to the first case. 
A very typical and most known example of
a selection theorem is the celebrated theorem of Michael. It
states that every lower semicontinuous (LSC) mapping $F:X \to Y$
from a paracompact domain $X$ into a Banach range space $Y$ admits
a continuous single-valued selection whenever each value $F(x), x
\in X,$ is a nonempty convex and closed subset of $Y.$

Consider now two multivalued mappings $F_{1}:X \to Y_{1}$,
$F_{2}:X \to Y_{2}$ and a single-valued mapping $L:Y_{1} \times
Y_{2} \to Y$. Denote by $L(F_{1};F_{2})$ the composite mapping,
which associates to each $x \in X$ the set
$$
\{y \in Y:y=L(y_{1};y_{2}),\,\,y_{1} \in F_{1}(x),\,\,y_{2} \in
F_{2}(x)\}.
$$

\proclaim{Definition 1.1} Let $f$ be a selection of the composite
mapping $L(F_{1};F_{2})$. A pair $(f_{1},f_{2})$ is said to be a
{\bf splitting} of $f$ if $f_{1}$ is a selection of $F_{1}$,
$f_{2}$ is a selection of $F_{2}$ and $f=L(f_{1};f_{2})$.
\endproclaim

In Sections 2 and 3 below we work in the category of topological
spaces. Thus the {\it splitting problem} (see \cite{RS1}) for the
triple $(F_{1}, F_{2}, L)$ is the problem of finding continuous
selections $f_{1}$ and $f_{2}$ which split a continuous selection
$f$ of the composite mapping $L(F_{1};F_{2})$.

For $Y_{1}=Y_{2}=Y$ and $L(y_{1};y_{2})=y_{1}+y_{2}$ we see the
specific problem of splitting into a sum of two items. More
generally, for constant multivalued mappings, the splitting problem
can be interpreted as the problem of continuous dependence of
solutions of the linear equation $y=L(y_{1};y_{2})$ on the data
$y$ and with constraints $y_{1} \in A$ and $y_{2}\in B$.

One more example: \,let $Y_{1}=Y_{2}=\R, F_{1}(\cdot)\equiv
[0,+\infty), F_{2}(\cdot)\equiv (-\infty; 0]$ and again
$L(y_{1};y_{2})=y_{1}+y_{2}$. Then $L(F_{1};F_{2})(\cdot) \equiv
\R$ and an arbitrary selection of $L(F_{1};F_{2})$ is simply an
arbitrary mapping from the
domain into $\R$. So in this case the solvability of
the splitting problem  means the existence of a
decomposition $f=f^{+}+f^{-}$, e.g. in the theory of the
Lebesgue integral\,(see \cite{H, Sect.25}).

Within the framework of the general theory of continuous
selections and due to the Banach open mapping principle it is
quite natural to restrict ourselves to the case of paracompact
domains $X$, Banach range spaces $Y_{1}, Y_{2}, Y$ and LSC
convex-valued and closed-valued mappings $F_{1},F_{2}$, and to the
case of linear continuous surjections $L:Y_{1} \times Y_{2} \to
Y$.

For a special case of the constant mappings $F_{1}(\cdot) \equiv A$
and $F_{1}(\cdot) \equiv B$, the splitting problem can be reduced
(Theorem 3.1)  to invertibility of a mapping $L:Y_{1} \times
Y_{2} \to Y$ with respect to an
appropriate family $\CC$ of subsets of
$Y_{1} \times Y_{2}$.

\proclaim{Definition 1.2} A linear continuous mapping $L:Z \to Y$
between Banach spaces is said to be $\CC$-hereditary invertible
for a family $\CC$ of subsets of $Z$ if for every $C \in \CC$ the
restriction $L|_{C}:C \to L(C)$ admits a right-inverse continuous
mapping $s: L(C) \to C$,\,\,$L|_{C}\circ s = id|_{L(C)}$.
\endproclaim

In terms of continuous selections,  $L:Z \to Y$ is
$\CC$-hereditary invertible whenever the inverse multivalued
mapping $(L|_{C})^{-1}: L(C) \to C$ admits a continuous selection.
Clearly, for a class $\CC$ consisting of closed and convex sets
the $\CC$-hereditary invertibility of $L:Z \to Y$ follows from
$\CC$- hereditary openess of $L$. This simply means that each
restriction $L|_{C}:C \to L(C)$ is an open mapping. Therefore
$\CC$-hereditary openess of $L$ guarantees that the Michael
selection theorem mentioned above is applicable to each mapping
$(L|_{C})^{-1}: L(C) \to C, \,\,C \in \CC$.

Unfortunately, as a rule $\CC$-hereditary openess (and also $\CC$-
hereditary invertibility) of an arbitrary map $L:Z \to Y$  is a
very restrictive property. For example, for the class $\CC$ of all
convex compacta this
means  that $\dim Z \leq 2$ or $\dim Y=1$ (Theorem
2.1 and Remark (1)). In Theorem 2.3 we prove that if the boundary
of a convex finite-dimensional compactum $C$ is "transversal" to
$\Ker L$ then $L|_{C}:C \to L(C)$ is an open mapping. On other hand,
finite dimensionality is here the principal point. Namely, Theorem
2.4 shows that in any infinite-dimensional Banach space $Z$ there
is a subcompactum $C$ for which all assumptions of Theorem 2.3.
hold, but $L|_{C}:C \to L(C)$ is not open, and moreover the
inverse mapping $(L|_{C})^{-1}: L(C) \to C$ admits no (even local)
continuous selection.

In Section 3 we apply positive results of Section 2 to finding of
the splittings. In particular, for a single-valued mapping $f$ to
a compact space $L(A,B)$ we obtain results on splitting of $f$
into mappings to $A$ and to $B$ (Theorems 3.5 and 3.6). As a
corollary, we prove that for the Minkowski sum $A+B$ of finite-dimensional
strictly convex bounded $A$ and $B$ the equality
$c=a(c)+b(c), c \in A+B$ holds for some continuous single-valued
mappings $a:A+B \to A$ and $b:A+B \to B$.

Finally, recall that the lower semicontinuity of a multivalued
mapping $F: X \to Y$ between topological spaces $X$ and $Y$ means
that for each points $x \in X$ and $y \in F(x),$ and each open
neighborhood $U(y)$, there exists an open neighborhood $V(x)$ such
that $F(x') \cap U(y) \ne \emptyset$, whenever $x' \in V(x)$. If
one identifies the mapping $F: X \to Y$ with its graph $\Gr F \ss X
\times Y$, then the lower semicontinuity of $F$ is equivalent to
the openess of the restriction $p_{1}|_{\Gr F}: \Gr F \to X$, where
$p_{1}:X \times Y \to X$ is the projector onto the first
coordinate. Roughly speaking, lower semicontinuous multivalued
mappings are exactly  inverses of open single-valued mappings.

\medskip
\medskip
\head {\bf 2. Hereditary openess and
invertibility}\endhead

\proclaim{Theorem 2.1} For any Banach space $Y$ the following statements are
equivalent:
\itemitem{(a)} each linear continuous surjection $L:Z \to Y$ from
a Banach space $Z$ is $\FF_{c}(Z)$-hereditary invertible with respect to
the family $\FF_{c}(Z)$ of all convex closed subsets of $Z$;
\itemitem{(b)} each linear continuous surjection $L:Z \to Y$ from
a Banach space $Z$ is $\CC_{c}$-hereditary invertible with respect to the
family $\CC_{c}$ of all convex subcompacta of $Z$; and
\itemitem{(c)} $\dim Y=1$.
\endproclaim
\demo{Proof} The implication $(a)\Rightarrow (b)$ is trivial. To
check $(b)\Rightarrow (c)$ we shall need the following
lemma.

\proclaim{Lemma 2.2}In the Euclidean 3-space $\R^{3}=\R^{2} \oplus \R^1$
there is a convex compact set $C$ such that the restriction $P|_{C}: C
\to P(C)$ of the orthogonal projection $P: \R^{3} \to \R^{2}$ is
not an open mapping. Moreover, the inverse multivalued mapping
$(P|_{C})^{-1}$ admits no continuous selection.
\endproclaim
\demo{Proof} Let $K$ be one full rotation
of the
spiral
$$
K=\{(\cos t, \sin t,t):0\leq t \leq 2\pi\}
$$
and $C=\conv  K=\overline{\conv }K$.  Suppose to the contrary that the
point-preimages multivalued mapping $(P|_{C})^{-1}: P(C) \to C$
admits a continuous selection, say $s: P(C) \to C$. Observe that
$P(C)$ is the unit disk $D=\{(r \cos t,r \sin t,0):0\leq t \leq
2\pi, 0 \leq r \leq 1 \}$ and that the mapping $(P|_{C})^{-1}$ is
single-valued over the whole boundary $\partial D$ of $D$ except
over the initial point $(1,0,0)$. Hence the continuous selection
$s$ coincides with $(P|_{C})^{-1}$ on $\partial D \backslash
\{(1,0,0)\}$.

Therefore $\lim_{t \to 0+0} s(\cos t,\sin t,0) = (1,0,0)$ and
$\lim_{t \to 2\pi-0} s(\cos t,\sin t,0) = (1,0,2\pi)$, which
contradicts the continuity of $s$. Note that in fact, the mapping
$(P|_{C})^{-1}$ admits no selections which is continuous at the
point $(1,0,0)$. \q\enddemo

Now, suppose that the assumption $(c)$ does not hold, i.e. $\dim Y \geq 2$.
Hence $Y=\R^{2} \oplus Y'$ for some Banach space $Y'$. Let
$Z=\R^{3} \oplus Y'$. Then we can map $\R^{3} \to \R^{2}$ as in
Lemma 2.2, and map $Y'$ onto $Y'$ identically, consider the direct sum of
these linear surjections and obtain a contradiction with the
assumption $(b)$ on the existence of the right inverse for the
restriction $P|_{C}: C \to P(C)$.

In order to prove $(c)\Rightarrow (a)$, let us first check that the
restriction $L|_{C}:C \to L(C)$ is an open mapping for every linear
continuous map $L:Z \to Y$ and for every convex set
$C \ss Z$. Suppose to the
contrary
that $L|_{C}$ is not open at some point $z \in C$. Then there
exist a number $\e>0$ and a sequence
$\{y_{n}\}_{n=1}^{\i}$ with $y_{n} \in L(C)$ such that $y_{n} \to
L(z), n \to \i$ and
$$
\dist(z; L^{-1}(y_{n}) \cap C) \geq  \e, \qquad n \in \N.
$$

The set $L(C)$ is convex and one-dimensional. Thus one can assume
that $\{y_{n}\}_{n=1}^{\i}$ is monotone. Let\,\,$y_{n} \to
L(z)+0$. Then $y_{n}=(1-t_{n})L(z)+t_{n}y_{1},  t_{n} \to 0+0$. By
the choice of $\{y_{n}\}_{n=1}^{\i}$ there exists a point $z_{1}
\in L^{-1}(y_{1}) \cap C$. Hence $z_{n}=(1-t_{n})z+t_{n}z_{1} \in
[z,z_{1}] \ss C$ and $L(z_{n})=y_{n}$. So $z_{n} \in L^{-1}(y_{n})
\cap C$ and $z_{n} \to z$. Thus $\dist(z; L^{-1}(y_{n}) \cap C) \to
0$ which contradicts the fact that $\dist(z; L^{-1}(y_{n}) \cap C)
\geq \e$.

Now, let us return to the case when $C \in \FF_{c}(Z)$. Since the
set $L(C)$ is metrizable and hence paracompact, all values of the
mapping $y \mapsto L^{-1}(y) \cap C$ are nonempty, convex and
closed. Such the mapping is LSC because $L|_{C}:C \to L(C)$ is an
open mapping. So applying the Michael selection theorem we find
the continuous right inverse of $L|_{C}:C \to L(C)$. \q\enddemo

{\bf Remarks:}

{\bf (1)} In the same way one can prove that $\CC_{c}$-hereditary
invertibility characterizes Banach spaces $Z$ with $\dim Z \leq 2$.

{\bf (2)} The analog of Theorem 2.1 holds under substitution of
hereditary openess instead of hereditary invertibility even
without closedness assumption for convex subsets of $Z$ in $(a)$.
In fact, one can use instead of the example from Lemma 2.2 another
(widely known) example of the convex hull $C \ss \R^{3}$ of the
set $\{(\cos t, \sin t, 0): 0 \leq t \leq 2\pi\} \cup \{(1,0,1)\}$
and orthogonal projection $p: \R^{3} \to \R^{2},
p(x,y,z)=(x,y,0)$. Note that $(p|_{C})^{-1}$ here admits the
obvious (identical) continuous selection. This is the key
difference with the example from Lemma 2.2.

\medskip

Theorem 2.1 shows that separate and independent assumptions on
linear mapping $L$ and on a convex compact set $C \ss Z$ cannot give
an essential result. So some linking properties on $L$ and $C$
are needed.

Let us recall that for a convex subset $M$ of a Banach space $Z$
there are (at least) two approaches to the notion of its
relatively inner point. First, a point $m \in M$ is said to be
inner (in the {\it metric} sense) point of $M$ provided that for some
positive $\e$ the intersection $D(m;\e)\cap \aff(M)$ is subset of
$M$. Here and below $D(m;\e)$ denotes the open ball with radius
$\e$ centered at $m$. Second, a point $m \in M$ is said to be
inner (in  the {\it convex} sense) point of $M$ provided that for each
$x \in M, x\neq m$, there exists $y \in M$ such that $m \in
[x;y)$. Here, $[x;y)$ is the straight line semiinterval, i.e. the
segment $[x;y]$ without the end point $y$.

A great advantage of finite-dimensional convex sets is that for
them these approaches are equivalent (see \cite{W,2.3.6 and
2.6.10}). The Hilbert cube $Q$, lying in any Banach, or Frechet
space, has no inner (in the {\it metric} sense) points. But $Q$
certainly has inner (in the {\it convex} sense) points: they
constitute the so-called {\it pseudo-interior} of the Hilbert
cube. Note that each infinite-dimensional convex compact subset of
a Frechet space is homeomorphic to $Q$, due to the the Keller
theorem \cite{Mi}.

Below we shall use this equivalence without any special reference
and we shall denote by $\int(M)$ (resp.,\,$\partial(M)$)\,\,the set of all
inner (resp.,\,boundary) points of a finite-dimensional convex set
$M$. Observe that $\int(A \times B)=\int(A) \times \int(B)$.

\proclaim {Theorem 2.3} Let $L:X \to Y$ be a linear continuous
surjection between Banach spaces. Let $C \ss X$ be a convex
finite-dimensional bounded subset of $X$ such that the boundary
$\partial(C)$ contains no segments parallel to the kernel
$\Ker(L)$. Then the restriction $L|_{C}:C \to L(C)$ is an open
mapping.
\endproclaim
\demo{Proof}

1) Let $x \in \int(C)$. Then the conclusion follows from the Banach
open mapping principle, applied to the restriction $L|_{\aff(C)}$.

2) Let $x \in C \cap \partial(C), L(x)=y$, but suppose that
$L^{-1}(y)$ intersects $\int(C)$. So let $x_{0} \in L^{-1}(y)\cap
\int(C)$. It is a well known and fundamental fact that the whole
semiinterval $(x;x_{0}]$ lies in $\int(C)$ \cite{W,2.3.4}. Thus for
every $\e>0$ there is $x_{\e} \in L^{-1}(y)\cap \int(C) \cap
D(x;\e)$. Choose $\d>0$ such that $D(x_{\e};\d) \ss D(x;\e)$. Due
to the case (1) the image $L(D(x_{\e};\d) \cap C)$ contains some
neighborhood, say $V(y)$ of the point $y$ in $L(C)$. This is why
$$
V(y) \ss  L(D(x_{\e};\d) \cap C)  \ss L(D(x;\e) \cap C),
$$
i.e. the restriction $L|_{C}:C \to L(C)$ is open at the point $x$.

3) Let $x \in C \cap \partial(C), L(x)=y$ and $L^{-1}(y) \cap
\int(C) = \varnothing$. Hence $(L^{-1}(y) \cap C) \ss \partial(C)$.
If $x_{1}$ and $x_{2}$ are two distinct points in $L^{-1}(y) \cap
C$ then the segment $[x_{1},x_{2}]$ is parallel to $\Ker(L)$ and
lies in the boundary $\partial(C)$. This is a contradiction with
the assumption of the theorem. Thus $L^{-1}(y) \cap C = \{x\}$.

Suppose to the contrary, that the mapping $L|_{C}$ is not open at
$x$. This means that for some $\e>0$ and for some sequence $y_{n}
\to y, y_{n} \in L(C), n\to \infty$,\, each distance
$\dist(L^{-1}(y_{n})\cap C, x)$ is greater than or equal to $\e$.
For each $n \in \N$ pick $x_{n} \in L^{-1}(y_{n})\cap C$. Due to
precompactness of $C$ we can choose convergent subsequence, say
$x_{n_{k}} \to x_{0}, x_{0} \in Cl(C), k \to \infty$. Then
$y_{n_{k}} = L(x_{n_{k}}) \to L(x_{0})$ and $L(x_{0})=y$, or
$x_{0} \in L^{-1}(y) \cap Cl(C)$. But $L^{-1}(y) \cap \int(C) =
\varnothing$. Therefore $x_{0} \in L^{-1}(y) \cap \partial(C)$ and
$x_{0}=x$, due to the transversality type assumption that the
boundary $\partial(C)$ contains no segments parallel to the kernel
$\Ker(L)$. So
$$
\dist(x_{n_{k}}, x) \geq \dist(L^{-1}(y_{n_{k}})\cap C, x) \to
0,\,\,\, k\to\infty
$$
Contradiction. \q
\enddemo

The following theorem demonstrates that without the finite
dimensionality the
restriction local invertibility can fail
even at the
inner points.

\proclaim{Theorem 2.4} For every infinite dimensional Banach space
$Z$ and for every continuous linear projector $P: Z \to Z$ with
$\dim \Ker(P)=1$ there is an infinite-dimensional convex compact
subset $C \ss Z$ and an inner point $z \in C$ such that the
restriction $P|_{C}: C \to P(C)$ is not open at $z$. Moreover, the
inverse multivalued mapping $(P|_{C})^{-1}: P(C) \to C$ admits no
continuous selections over an arbitrary neighborhood of the point
$P(z)$.
\endproclaim
\demo{Proof} Choose any basic Schauder normalized sequence $e_{1},
e_{2},...,e_{n},...\,,\,\,\|e_{n}\|=1$ with $e_{1} \in \Ker P$ and
$e_{n} \in \Im  P,\,\,n>1$. So, $P(e_{1})=0, P(e_{n})=e_{n},\,\,n>1$.
Define
$$
C=\overline{\conv } \left\{ \frac{e_{n}}{n},
\,\,2e_{1}-\frac{e_{n}}{n} \right\}_{n=1}^{\i}, \qquad K=\conv 
\left\{\frac{e_{n}}{n},\,\,
2e_{1}-\frac{e_{n}}{n}\right\}_{n=1}^{\i}.
$$

The set $\{\frac{e_{n}}{n}, 2e_{1}-\frac{e_{n}}{n}\}_{n=1}^{\i}$
is precompact because it consists of two convergent sequences
$\frac{e_{n}}{n} \to 0$ and $2e_{1}-\frac{e_{n}}{n} \to 2e_{1}$.
Hence $K$ is also precompact set and $C$ is a convex compact subset
of $X$. The point $e_{1}$ is the center of symmetry of the set $C$
and hence is its inner point.

Suppose that we  have already checked that the multivalued mapping
$(P|_{C})^{-1}: P(C) \to C$ is single-valued over the set
$\{\frac{e_{n}}{n}, 2e_{1}-\frac{e_{n}}{n}\}_{n=1}^{\i}$. In other
words, suppose that we have proved that
$$
(P|_{C})^{-1}\left(\frac{e_{n}}{n}\right)=\frac{e_{n}}{n}, \qquad
(P|_{C})^{-1}\left(-\frac{e_{n}}{n}\right)=2e_{1}-\frac{e_{n}}{n}.
$$
In this assumption, if $s:U \cap P(C) \to C$ is a continuous
selection of $(P|_{C})^{-1}$ over some neighborhood $U$ of the
point $P(e_{1})=0$ then
$$
\lim_{n \to \i} s\left(\frac{e_{n}}{n}\right) = 0, \qquad \lim_{n
\to \i} s\left(-\frac{e_{n}}{n}\right)=2e_{1},
$$
which contradicts the continuity of $s$ at $0$.

In order to complete the proof it suffices to check that $\l
e_{1}+\frac{e_{n}}{n} \in C$ if and only if  $\l=0$ and, analogously $\mu
e_{1}+\left(2e_{1}-\frac{e_{n}}{n}\right) \in C$ if and only if  $\mu=0$.

\proclaim {Lemma 2.5} For each \,\,$z=\l
e_{1}+\frac{e_{n}}{n},\,\,\, \l>0, n>1$ there exists $d>0$ such
that\,\,\, $\dist(z,K)\geq d$ and hence \,\,\,$\dist(z,C)\geq d$.
\endproclaim
\demo{Proof} For every $N \in \N$,\, let
$K_{N}=\conv \left\{\frac{e_{k}}{k},
2e_{1}-\frac{e_{k}}{k}\right\}_{k=1}^{N}$. Then $K=\bigcup
_{n=1}^{N} K_{N}$ and $\dist(z,K)=\inf\{\dist(z,K_{N}): N \in \N
\}$. Clearly,\, $K_{n-1} \ss \span \{e_{k}\}_{k=1}^{n-1}$ and $z \notin
\span \{e_{k}\}_{k=1}^{n-1}$. Therefore $z \notin K_{n-1}$ and
$\dist(z,K_{n-1}) = d_{1}>0$.

Thus we must consider the case $N \geq n$. So let $y \in K_{N}$,
i.e.
$$
y=\a_{1}e_{1}+\sum_{k=2}^{N}
\left(\a_{k}\frac{e_{k}}{k}+\b_{k}(2e_{1}-\frac{e_{k}}{k})\right)=
\left(\a_{1}+2\sum_{k=2}^{N}
\b_{k}\right)e_{1}+\sum_{k=2}^{N}(\a_{k}-\b_{k})\frac{e_{k}}{k}
$$
for some nonnegative $\a_{1}, \a_{2},...\a_{N}, \b_{2},...,\b_{N}$
with $\a_{1}+\sum_{k=2}^{N}(\a_{k}+\b_{k})=1$.

Hence
$$
z-y\,\,=\,\,\left(\l-(\a_{1}+2\sum_{k=2}^{N}
\b_{k})\right)e_{1}+\left(1-(\a_{n}-\b_{n})\right)\frac{e_{n}}{n}+\sum_{k=2,\,k\neq
n}^{N}(\b_{k}-\a_{k})\frac{e_{k}}{k}.
$$
Note that $\a_{n}-\b_{n} \leq \a_{n} \leq 1$ and therefore
$1-(\a_{n}-\b_{n}) \geq 0$.

{\it Case (a).} Let $1-(\a_{n}-\b_{n}) \geq \f{\l}{3}$. Then
$$
\|z-y\| \quad \geq \quad \f{1-(\a_{n}-\b_{n})}{n\|P_{n}\|} \quad
\geq \quad \f{\l}{3n\|P_{n}\|},
$$
where $P_{n}: \overline{\span }\{e_{k}\}_{k=1}^{\i} \to
\span \{e_{n}\}$ is the continuous linear projection onto the $n-$th
coordinate. Recall that $\|P_{n}(u)\| \leq \|P_{n}\|\cdot \|u\|$.
In our case $u=z-y$ and $P_{n}(u) = (1-(\a_{n}-\b_{n}))e_{n}$.

{\it Case (b).}\,\,Let $1-(\a_{n}-\b_{n}) < \f{\l}{3}$. Then $1-
\f{\l}{3}< \a_{n}-\b_{n} \leq \a_{n}$ and $1-\a_{n}<\f{\l}{3}$.
So
$$
\b_{2}+...+\b_{N} \leq \a_{1}+\b_{2}+...+\b_{N} \leq
1-\a_{n}<\f{\l}{3}
$$
and $\a_{1}+2\sum_{k=2}^{N} \b_{k} < \f{2\l}{3}$. Hence $\l -
\left(\a_{1}+2\sum_{k=2}^{N} \b_{k}\right) > \f{\l}{3}$ and $
\|x-y\| > \f{\l}{3\|P_{1}\|},$ where $P_{1}:
\overline{\span }\{e_{k}\}_{k=1}^{\i} \to \span \{e_{1}\}$ is the
continuous linear projection onto the first coordinate.

Thus, in any case $ \|z-y\|  \geq \min\left\{\f{\l}{3n\|P_{n}\|},
\f{\l}{3\|P_{1}\|}\right\}=d_{2}>0.$ Finally, for each $N \in \N$
$$
\dist(z,K_{N})= \inf\{\|z-y\|: y \in K_{N} \} \geq \min\{d_{1},
d_{2}\} = d >0
$$
and $\dist(z,K) \geq d$. This completes the proof of the lemma and
also of the theorem. \q\enddemo
\enddemo

\head {\bf 3. Splitting selections}\endhead

To reduce the splitting problem for constant multivalued mappings to
the hereditary invertibility property the following simple
statement is useful. We state it in a rather abstract form.

\proclaim{Theorem 3.1} Suppose that a continuous surjection
$L:Y_{1} \times Y_{2} \to Y$ between Banach spaces is $(\CC_{1}
\times \CC_{2})$-hereditary invertible with respect to some families
$\CC_{1}$ subsets of $Y_{1}$ and $\CC_{2}$ subsets of $Y_{2}$. Let
$A \in \CC_{1}$ and $B \in \CC_{2}$.  Then the splitting problem
for the triple $(F_{1}(\cdot) \equiv A, F_{2}(\cdot) \equiv B, L)$
is solvable for an arbitrary domain $X$.
\endproclaim
\demo{Proof} Under the assumptions the composite mapping
$F=L(F_{1},F_{2})$ is the constant multivalued mapping
$F(\cdot)\equiv L(A,B)$. Hence its continuous single-valued
selection, say $f$, simply is an arbitrary continuous
single-valued mapping $f: X \to L(A,B)$ from a domain $X$. We
define an auxiliary multivalued mapping $\Phi:X\to Y_{1} \times
Y_{2}$, by setting
$$
\Phi(x)=\{(y_{1};y_{2})|\,\,y_{1} \in A, y_{2} \in B,
\,\,\,L(y_{1};y_{2})=f(x)\} = (A\times B) \cap L^{-1}(f(x)).
$$
All values of $\Phi$ are nonempty because $f(X) \ss L(A,B)$. So
first, to each argument $x \in X$  in the continuous fashion
there corresponds the point $y=f(x) \in L(A,B) \ss Y$. Second, to this
point $y$ corresponds the set $(A\times B) \cap L^{-1}(y)$. And
the $(\CC_{1} \times \CC_{2})$-hereditary invertibility of $L$
means exactly that the last multivalued
correspondence has a continuous selection (see Definition 1.2). Thus its composition
with $f$ is a continuous selection, say $\phi$, of $\Phi$.

So if $f_{1}=p_{1} \circ \phi: X \to A$ and $f_{2}=p_{2} \circ
\phi: X \to B$, where $p_{i}(y_{1},y_{2})=y_{i}, i=1,\,2$ are
"coordinate" projections $p_{i}: Y_{1} \times Y_{2} \to Y_{i}$,
then
$$
L(f_{1}(x), f_{2}(x))= L(p_{1}(\phi(x)),
p_{2}(\phi(x)))=L(\phi(x))=f(x),
$$
because $\phi(x) \in \Phi(x) \ss L^{-1}(f(x)), \,\,\, x \in X.$
Thus the pair $(f_{1}, f_{2})$ splits the mapping $f$.  \q\enddemo

We emphasize that Theorem 3.1 is a conditional statement which
simply reduces one problem to another: the checking of $(\CC_{1}
\times \CC_{2})$-hereditary invertibility is a separate and
nontrivial job. Theorem 3.1 gives a way of transferring the
results from the previous section to splitting of continuous
selections. First we transfer the example from Lemma 2.2.

\proclaim{Example 3.2} For any 2-dimensional cell $D$ there exist:
\itemitem{(a)} constant multivalued mappings $F_{1}:D \to \R^{3}$
and  $F_{2}:D \to \R$ with convex compact values; \itemitem{(b)} a
linear surjection $L: \R^{3} \oplus \R \to \R^{2}$; and
\itemitem{(c)} a continuous selection $f$ of the composite mapping
$F=L(F_{1},F_{2})$, such that $f\neq L(f_{1},f_{2})$ for any
continuous selections $f_{i}$ of $F_{i}$, i=1,2.
\endproclaim
\demo{Proof} In the notations of Lemma 2.2 let $$D=P(C),
F_{1}(\cdot)\equiv C, F_{2}(\cdot) \equiv [0;1],   L=P\oplus
0|_{\R}:\R^{3} \oplus \R \to \R^{2}.$$

Then $L(C\oplus [0;1]) = D$, $F(\cdot)=L(F_{1},F_{2})(\cdot)
\equiv D$ and $f=id|_{D}$ is a continuous selection of $F$.
Suppose to the contrary that $f=L(f_{1},f_{2})$ for some
continuous selections $f_{i}$ of $F_{i}$, i.e. for mappings
$f_{1}:D \to C$ and $f_{2}:D \to [0;1]$. But the surjection $L$
"forgets" the second coordinate. Hence
$$
D=f(x)=L(f_{1}(x),f_{2}(x))=P(f_{1}(x)), x \in D
$$
or $f_{1}(x) \in C \cap P^{-1}(x)$.

This means that $f_{1}$ is a continuous selection of multivalued
mapping $x\mapsto C \cap P^{-1}(x), x \in D$ which contradicts
Lemma 2.2.\q\enddemo

For application of Theorem 2.3 we need some additional
smoothness-like restriction on boundaries of convex sets (compare
with the notion of a strictly convex Banach space).

\proclaim{Definition 3.3} The convex subset $C$ of a Banach space
is said to be {\bf strictly convex} if the middle point of any
nontrivial segment $[x,y],\,\,x \in C, y \in C$ is an inner (in the
convex sense) point of $C$.
\endproclaim

Equivalently, the boundary of $C$ contains no straight line
segment.

\proclaim {Theorem 3.4} Let $A$ and $B$ be strictly convex
finite-dimensional bounded subsets of  Banach spaces $Y_{1}$ and
$Y_{2}$, respectively. Let $L:Y_{1} \times Y_{2} \to Y$ be a
linear continuous surjection with kernel $\Ker(L)$ transversal to
$Y_{1} \times \{0\}$ and $\{0\} \times Y_{2} $. Then the
restriction $L|_{A \times B}:A \times B \to L(A \times B)$ is an
open mapping.
\endproclaim
\demo{Proof} In view of Theorem 2.3, it suffices to check only
that the boundary $\partial(A \times B)$ contains no segment
parallel to $\Ker(L)$. Suppose to the contrary that $c_{1}\neq
c_{2}$, $[c_{1},c_{2}]=[(a_{1}, b_{1}),(a_{2}, b_{2})] \ss
\partial(A \times B)$ and $[c_{1},c_{2}]$ is parallel to $\Ker(L)$.
This means that $(a_{1}-a_{2}, b_{1}-b_{2}) \in \Ker(L)$. So if
$a_{1}=a_{2}$ then the transversality assumption implies that
$b_{1}=b_{2}$ and hence $c_{1}=c_{2}$. Contradiction. Hence
$a_{1}\neq a_{2}$ and analogously $b_{1}\neq b_{2}$.

%Next, equality $L(a_{1}, b_{1})=L(a_{2}, b_{2})$ and linearity of
%$L$ show that $L(a_{1}, b_{1})=L(a_{2}, b_{2})=
%L(0,5(a_{1}+a_{2}), 0,5(b_{1}+b_{2}))$.

By strict convexity $a'=0,5(a_{1}+a_{2}) \in \int(A)$ and
$b'=0,5(b_{1}+b_{2}) \in \int(B)$. But $(a',b') \in [c_{1},c_{2}]$.
So the segment $[c_{1},c_{2}]$ intersects $\int(A \times B)$ which
contradicts the existence of inclusion $[c_{1},c_{2}] \ss \partial(A \times
B)$. \q \enddemo

Theorems 3.1 and 3.4 together imply:

\proclaim {Theorem 3.5} Let $A$ and $B$ be strictly convex
finite-dimensional bounded subsets of Banach spaces $Y_{1}$ and
$Y_{2}$, respectively. Let $L:Y_{1} \times Y_{2} \to Y$ be a
linear continuous surjection with kernel $\Ker(L)$ transversal to
$Y_{1} \times \{0\}$ and $\{0\} \times Y_{2} $. Then for every
continuous single-valued mapping $f: X \to L(A,B)$ from a domain
$X$ there are continuous single-valued mappings $f_{1}: X \to A$
and $f_{2}: X \to B$ such that
$$
L(f_{1}(x), f_{2}(x))=f(x), \qquad x \in X.
$$
\endproclaim
\demo{Proof} Theorem 3.4 implies that the restriction $L|_{A
\times B}: A \times B \to L(A \times B)$ is an open mapping. Its
image is a metric (and hence, perfectly normal) space. All its
point-preimages are nonempty convex finite-dimensional subsets of
a separable (finite-dimensional, in fact)Banach space $\span (A)
\times \span (B)$. Hence Theorem 3.1''' from \cite{M} shows that
$L|_{A \times B}: A \times B \to L(A \times B)$ is
$\CC$-hereditary invertible, where $\CC$ is the family of all
strictly convex finite-dimensional subsets of a Banach space
$\span (A) \times \span (B)$. So an application of Theorem 3.1
completes the proof. \q
\enddemo

Remark that for a convex closed-valued LSC mappings $F_{1}:X \to
Y_{1}$ and $F_{2}:X \to Y_{2}$ and for a linear continuous
surjection $L:Y_{1} \times Y_{2} \to Y$, the splitting problem has
an affirmative solution in the case of one-dimensional $Y_{1}$ and
$Y_{2}$ and arbitrary paracompact domains (see \cite{RS1, Theorem
3.1}). But in general, splittings of continuous selections exist
only if members of the triple $(F_{1}, F_{2}, L)$ properly agree.
See \cite{RS1, Example 4.2} for a counterexample even for the case
$\dim Y_{1}=2, \dim Y_{2}=1$ and for a countable domain.

We conclude the section by showing the partial case of Theorem 3.5
applying it for the Minkowski sum of convex sets.

\proclaim {Corollary 3.6} Let $A$ and $B$ be strictly convex
finite-dimensional bounded subsets of  Banach spaces $Y$. Then
there are continuous single-valued mappings $a:A+B \to A$ and
$b:A+B \to B$ such that $c=a(c)+b(c)$ for all $c \in A+B$.
\endproclaim
\demo{Proof} In assumptions of Theorem 3.5 we choose the very
special linear continuous surjection $L:Y_{1} \times Y_{2} \to Y$
and special perfectly normal (in fact, metric) domain $C$. Namely,
$Y_{1}=Y_{2}=Y$, $L(y_{1},y_{2})=y_{1}+y_{2}$ and $C=A+B$.

Clearly $(y_{1},0) \in \Ker(L) \Leftrightarrow y_{1}=0$, i.e.\,the
kernel $\Ker(L)$ is transversal to $Y \times \{0\}$ and to $\{0\}
\times Y $. So Theorem 3.6 implies that the identity mapping $id:
C \to C$ admits a splitting $id=L(f_{1},f_{2})$ for some
continuous single-valued $a:C \to A$ and $b:C \to B$.

In other words, if $c \in C$ and $c \mapsto \{(a,b): c=a+b\}$ then
we can always assume that $a=a(c)$ and $b=b(c)$ are continuous
items with respect to the data $c \in C$. \q
\enddemo

Analogously, the another version of Theorem 3.5 states that the
continuous mapping $f$ from $X$ to the Minkowski sum $A+B$ splits
into a sum of two continuous mappings $f_{1}:X \to A$ and $f_{2}:X
\to B$, whenever $A$ and $B$ are strictly convex
finite-dimensional bounded subsets of a Banach spaces $Y$.

Finally, we guess that the strict convexity assumption can be
weakened in some ways, but that in general, Corollary 3.6 does
not hold for an arbitrary convex finite-dimensional compacta.

\head {\bf Acknowledgements} \endhead
The first author was supported by the Slovenian Research Agency
grants No. P1-0292-0101-04 and Bl-RU/05-07/7. The second author
was supported by the RFBR grant No.\,05-01-00993.
We thank the referee for several comments and suggestions.


\medskip

\Refs
\widestnumber\key{10000}

\ref \key{AF} \by J.-P. Aubin and H. Frankowska \book Set-Valued
Analysis \publ Birkh\"{a}user \publaddr Basel \yr 1990
\endref

\ref \key{H} \by P. Halmos \book Measure Theory \publ Van Nostrand
Reinhold Company \publaddr New York \yr 1969
\endref

\ref \key{M} \by E. Michael \paper Continuous selections. I \jour
Ann. of Math. \vol{\bf 63} \yr 1956 \pages 361--382
\endref

\ref \key{M3} \by E. Michael \paper Paraconvex sets \jour
Scand.\,Math. \vol{\bf 7} \yr 1959 \pages 372--376 \endref

\ref \key{Mi} \by J. van Mill \book Infinite-Dimensional Topology.
Prerequisites and Introduction \publ North Holland Company
\publaddr Amsterdam \yr 1989
\endref


\ref \key{RS} \by D. Repov\v{s} and P. V. Semenov \book Continuous
Selections of Multivalued Mappings, {\rm Mathematics and Its
Applications} \vol {\bf 455} \publ Kluwer \publaddr Dordrecht \yr
1998 \endref

\ref \key{RS1} \by D. Repov\v{s} and P. V. Semenov \paper Sections
of convex bodies and splitting problem for selections \jour J.
Math. Anal. Appl., !!!ADD THE DATA!!!
\endref

\ref \key{W} \by R. Webster \book Convexity \publ Oxford Science
 Publ\publaddr New York \yr 1994
\endref


\endRefs
\enddocument